\documentclass[10pt,twoside,leqno,twocolumn]{article}
\usepackage{ltexpprt}
\usepackage{graphicx}
\usepackage{amsmath}
\usepackage{balance}
\usepackage{hyperref}
\usepackage{mathtools}
\usepackage{bbm}
\usepackage[ruled,vlined]{algorithm2e}
\usepackage[usenames,dvipsnames]{color}
\usepackage[font={small}]{caption}
\usepackage{geometry}

\newcommand{\newMargin}{0.86in}
\newgeometry{left=\newMargin,right=\newMargin,top=1in,bottom=1in}

\usepackage[hyphenbreaks]{breakurl}

\AtBeginDocument{
	\abovedisplayskip=1.5ex plus 1.5ex minus 1.5ex
	\belowdisplayskip=1.5ex plus 1.5ex minus 1.5ex
	\abovedisplayshortskip=1.5ex plus 0.5ex minus 1ex
	\belowdisplayshortskip=1.5ex plus 0.5ex minus 1ex
}

\makeatletter
\def\footnoterule{
  \hrule \@width 1.0cm \kern 1.6\p@}
\makeatother

\graphicspath{{fig/}{dat/}}

\makeatletter
\providecommand*{\input@path}{}
\g@addto@macro\input@path{{tex/}{tab/}{dat/}{fig/}{alg/}}
\makeatother

\allowdisplaybreaks
\newcommand{\ratio}[2]{#1\;:\:#2}
\newcommand{\about}[1]{|_{#1}}

\newcommand{\step}{\alpha}

\newcommand{\x}{\mathbf{x}_i}

\newcommand{\rr}{\mathbf{r}}

\newcommand{\w}{\mathbf{w}}
\newcommand{\wk}{\w_{k}}

\newcommand{\loss}{\mathcal{L} }

\newcommand{\p}{\mathbf{p} }

\newcommand{\polyDegree}{d}
\newcommand{\polyC}{c}
\newcommand{\polyCoeffScalar}{\polyC_\ell}
\newcommand{\polyCoeff}{\mathbf{\polyC} }
\newcommand{\errPoly}{\varepsilon_{d+1}}

\newcommand{\tol}{\theta}
\newcommand{\err}{\epsilon}

\newcommand{\W}{W}

\newcommand{\dphi}{\phi^{\prime}}

\newcommand{\Set}[1]{{\left\{#1\right\}}}

\newcommand{\abs}[1]{\left\lvert #1 \right\rvert}

\newcommand{\norm}[1]{\|#1\|}
\newcommand{\nnz}[1]{\text{nnz}\left(#1\right) }

\newcommand{\kdda}{{\sc kdd-a}}
\newcommand{\kddb}{{\sc kdd-b}}
\newcommand{\urldata}{{\sc url}}
\newcommand{\rcv}{{\sc rcv\footnotesize{1}}}
\newcommand{\eps}{{\sc e}psilon}

\newcommand{\costTime}{\tau_{\text{fg}}}
\newcommand{\iterTime}{\tau_{\text{tot}}}
\newcommand{\polyEvalTime}{\tau_{\ell\text{s}}}
\newcommand{\numEvals}{n_e}

\newcommand{\Ops}{\bigoplus}
\newcommand{\nnodes}{n_p}

\newcommand{\Real}{\mathbbm{R}}

\newcommand{\eqn}[1]{(\ref{#1})}

\begin{document}

	\author{
		Michael B Hynes\thanks{
			Department of Applied Mathematics,
			University of Waterloo,
			Waterloo, Canada.
			Email: mbhynes@uwaterloo.ca
		}
		\and
		Hans De Sterck\thanks{
			School of Mathematical Sciences,
			Monash University,
			Melbourne, Australia.
			Email: hans.desterck@monash.edu
		}
	}
	\date{}
	\title{ 
		\vspace{-1.20cm}
		\Large
		A polynomial expansion line search for large-scale 
		unconstrained minimization of smooth
		$L_2$-regularized loss functions, with
		implementation in Apache Spark\thanks{
			Supported in part by the National Sciences and
			Engineering Research Council through a Canada Graduate
			Scholarship
					}
	}
	\maketitle

\begin{abstract}

{
\bf
\footnotesize
In large-scale unconstrained optimization algorithms such as
limited memory BFGS (LBFGS), a common subproblem is a line
search minimizing the loss function along a descent direction.
Commonly used line searches iteratively find an
approximate solution for which the Wolfe conditions are
satisfied, typically requiring multiple function and
gradient evaluations per line search, which is expensive in
parallel due to communication requirements.
In this paper we propose a new line search
approach for cases where the loss
function is analytic, as in least squares regression,
logistic regression, or low rank matrix factorization.
We approximate the loss function by a truncated Taylor
polynomial, whose coefficients may be computed efficiently
in parallel with less communication than evaluating the
gradient, after which this polynomial may be minimized
with high accuracy in a neighbourhood of the expansion point.
The expansion may be repeated iteratively in a line search
invocation until the expansion point and minimum are sufficiently accurate.
Our Polynomial Expansion Line Search (PELS) was implemented in the
Apache Spark framework and used to accelerate the training of a
logistic regression model on binary classification datasets from the
LIBSVM repository
with LBFGS\ and the Nonlinear Conjugate Gradient (NCG)
method.
In large-scale numerical experiments in parallel on a
16-node cluster with 256 cores using the {\sc url}, {\sc kdd-a},
and {\sc kdd-b}\ datasets, the PELS\ approach
produced significant convergence improvements compared to the use of
classical Wolfe approximate line searches.
For example, to reach the final training label prediction
accuracies, LBFGS\ using PELS\ had speedup factors of
1.8--2 over LBFGS\ using a Wolfe approximate line search,
measured by both
the number of
iterations and the time required, due to the better accuracy
of step sizes computed in the line search.
PELS\ has the potential to significantly
accelerate
widely-used
parallel large-scale regression and factorization
computations, and is applicable
to important classes  of continuous optimization problems with
smooth loss functions.

}
\end{abstract}
\section{Introduction.}
\label{sec:intro}
In large-scale optimization and machine learning, a common
problem is the minimization of a loss
function $\mathcal{L} (\mathbf{w})$ with parameter vector $\mathbf{w} \in \mathbbm{R}^m$
\cite{hastie2001elements}.
For a set of $n$ observations $\mathcal{R} = {\left\{(\mathbf{x}_i,y_i)\right\}}_{i =
1}^{n}$ with $\mathbf{x}_i \in \mathbbm{R}^{p}$ and $y_i \in \mathbbm{R}$,
$\mathcal{L} (\mathbf{w})$ is expressed as the mean of the individual
losses ${\left\{f(\mathbf{w};\mathbf{x}_i,y_i) \right\}}$ evaluated at each observation with an
additional regularization term $\lambda R (\mathbf{w})$:
\begin{equation}\label{loss}
	\mathcal{L} (\mathbf{w}) = \lambda
	R (\mathbf{w}) + \frac{1}{n}\sum_{i=1}^n f(\mathbf{w};\mathbf{x}_i,y_i) .
\end{equation}
Here, $\lambda$ is a fitting parameter 
that tunes the magnitude of
the regularization penalty relative to the
losses \cite{yuan2012recent}.
Often the regularization is an $L_2$ penalty,
$\frac{1}{2} \|\mathbf{w}\|^2_2$, such that $R (\mathbf{w})$ is
smooth with respect to $\mathbf{w}$. 
The optimization problem for minimizing (\ref{loss}) is
\begin{equation}\label{lossOpt}
	\mathbf{w}^* = \arg\min_{\mathbf{w}} \mathcal{L} (\mathbf{w}), 
\end{equation} 
for which any unconstrained optimization algorithm may be
applied.
The use of nonsmooth regularization such as $L_1$
penalties that induce sparsity is also an active research topic,
however $L_2$ penalties are commonly used 
for many practical problems \cite{yuan2012recent,yuan2010comparison}
and improving algorithmic
performance in this case remains of significant interest.

The present work is situated in the context of batch or
minibatch optimization algorithms that solve
(\ref{lossOpt}) through iterative
updates with line searches, such as Gradient Descent (GD),
LBFGS, truncated Newton algorithms, or the Nonlinear
Conjugate Gradient (NCG) method (for a description of
these methods, see e.g.~Nocedal and Wright
\cite{nocedal2006numerical}).
These algorithms have an update for the
approximate solution $\mathbf{w}_{k}$ in the $k$th iteration as
\begin{equation}\label{update}
	\mathbf{w}_{k+1} = \mathbf{w}_{k} + \alpha_k\mathbf{p} _k 
\end{equation}
where $\mathbf{p} _k  \in \mathbbm{R}^m$ is the search direction and
the scalar $\alpha_k$ is a \emph{step size} such that
$\mathcal{L} (\mathbf{w}_{k+1})$ is either exactly or approximately minimized
in a univariate line search
along the ray $\mathbf{w}_{k} + \alpha \mathbf{p} _k $ with $\alpha > 0$.
Approximate
solutions that satisfy the Wolfe conditions are sought
in practice at the expense of poorer
convergence, since an inexact line search cannot in general be expected
to perform better than an exact line
search \cite{nocedal2006numerical}.
Since many commonly used univariate line search algorithms
require evaluating both $\mathcal{L} $ and $\mathbf{\nabla}\mathcal{L} $ in each line search
iteration, they
can be an expensive
component of batch optimization algorithms.

Our main idea in this paper is simple but powerful: when
(\ref{loss}) is smooth (infinitely differentiable), we
propose to approximate $\mathcal{L} (\mathbf{w} + \alpha\mathbf{p} )$ in a univariate
line search by a low-degree polynomial expansion, and show
that the coefficients of this polynomial may be computed in
a single pass over the dataset with modest
communication requirements, after which the
polynomial approximation may be minimized with high accuracy.
In each line search invocation, the expansion may be
repeated iteratively until the expansion point and
minimum are sufficiently accurate.
The advantages of our Polynomial Expansion Line Search\ (PELS) stem from
two main improvements.
Firstly, the PELS\ technique
can obtain more accurate minima when the high intrinsic potential accuracy of the
polynomial expansion is realized, which may lead to significantly
fewer iterations of the
optimization method to reach a desired accuracy than
with classical approximate line searches. Secondly, if multiple iterations are
required in a line search, the PELS\ method is much more efficient in terms of
parallel communication than the iterations in
approximate line searches that seek to impose the Wolfe conditions and require
evaluating $\mathbf{\nabla}\mathcal{L} $ in each line search iteration:
aggregating the
polynomial coefficients requires much less communication than
aggregating the $m$-dimensional gradient vectors
${\left\{\mathbf{\nabla} f(\mathbf{w};\mathbf{x}_i,y_i)\right\}}$.
Additionally, when (\ref{loss}) is polynomial,
the PELS\ expansion can be made \emph{exact} with a 
sufficiently large degree.

The contributions of this paper are organized as follows.
In \S\ref{sec:background}, the PELS\ algorithm is presented for
smooth loss functions, 
and a detailed example is given for $L_2$-regularized logistic
regression.
In \S\ref{sec:methods} we detail an implementation in the Apache Spark framework
\cite{zaharia2012resilient} for fault-tolerant distributed
computing, where it has been used to accelerate the training
of logistic regression models by GD, NCG, and
LBFGS\ on several large binary classification datasets.
Finally, \S\ref{sec:results} compares the performance of
these algorithms in Apache Spark using both
PELS\
and a standard approximate line search
scheme that may use multiple evaluations of $\mathcal{L} $
and $\mathbf{\nabla}\mathcal{L} $ per line search invocation.

\section{Background \& Theory.}
\label{sec:background}

\subsection{Line Search Formulation.}
\label{sec:ls}
A univariate line search is the unconstrained minimization problem of
determining a step size $\alpha^*$ that minimizes $\mathcal{L} (\mathbf{w})$
along a fixed descent direction $\mathbf{p} $, formulated as
\begin{equation}\label{ls}
	\alpha^* = \arg\min_{\alpha > 0} \mathcal{L} (\mathbf{w} + \alpha\mathbf{p} ) =
	\arg\min_{\alpha>0}\phi(\alpha).
\end{equation}
Since computing $\mathcal{L} (\mathbf{w} + \alpha\mathbf{p} )$ in (\ref{ls}) is an expensive
operation that requires evaluating 
$f(\mathbf{w} + \alpha\mathbf{p} ;\mathbf{x}_i,y_i) $ for each observation
$(\mathbf{x}_i,y_i)$, an approximate solution to (\ref{ls}) is often
sought instead of an exact solution in order to reduce the number of
function evaluations required in the line search.
Approximate line searches compute a sequence of iterates
${\left\{\alpha_j\right\}}_{j \geq 0}$ until convergence criteria are
satisfied, which are normally
the Strong Wolfe conditions \cite{nocedal2006numerical}.
These conditions are
a set of two
inequalities guaranteeing (i) sufficient decrease as
\begin{equation}\label{decrease}
	\phi(\alpha) \leq \phi(0) +
	\nu_1\alpha\phi^{\prime}(0),
\end{equation}
and (ii) a curvature condition requiring
\begin{equation}\label{curvature}
	\left\lvert \phi^{\prime}(\alpha) \right\rvert \leq \nu_2
	\left\lvert \phi^{\prime}(0) \right\rvert
\end{equation}
for $0<\nu_1<\nu_2<1$ (where $\nu_1 \approx 10^{-4}$ and
$\nu_2
\approx 0.9$ \cite{nocedal2006numerical}). 

There are a multitude of widely-used inexact univariate line search
algorithms for satisfying (\ref{decrease}) and
(\ref{curvature})
for general $\mathcal{L} (\mathbf{w})$ 
(see e.g.~\cite{more1994line,al1986efficient,hager1989derivative,hager2005new}), 
and we refer to an algorithm in this class as a
\emph{Wolfe approximate} (WA) line
search.
Most successful WA\ algorithms are variants of the
following 
classic interpolation scheme, summarized in Alg.\;\ref{alg:interpLS}.
In each $j$th iteration of the line search, an interval $I_j =
[\alpha_l,\alpha_u]$ containing $\alpha^*$ is determined, and the next
$\alpha_{j+1}$ is generated by the minimization of an
interpolated cubic polynomial $P_j(\alpha)$ with control points 
$\phi(\alpha_l),\phi(\alpha_u),\phi^{\prime}(\alpha_l)$,
and $\phi^{\prime}(\alpha_u)$ (i.e.~$P_j(\alpha_l) =
\phi(\alpha_l)$, $P_j^{\prime}(\alpha_l) =
\phi^{\prime}(\alpha_l)$, etc).
The next interval $I_{j+1}$ is computed such that the
endpoints are closer to $\alpha^*$, and methods such as
bisection, secant, and dual interpolation/minimization
can be used to shrink the interval.
Each evaluation of $\phi(\alpha)$ and $\phi^{\prime}(\alpha)$
requires an $\mathcal{O}\left[ n \right]$ pass over $\mathcal{R}$, and it is
typical in publicly available codes to structure
optimization routines with a function that computes both
$\phi(\alpha)$ and $\phi^{\prime}(\alpha)$ simultaneously
by evaluating $\mathcal{L} (\mathbf{w} + \alpha\mathbf{p} )$ and $\mathbf{\nabla}\mathcal{L} (\mathbf{w} + \alpha\mathbf{p} )$, 
explicitly returning a scalar and an $m$-dimensional
vector such that, if the current step size is accepted, the
computed value of $\mathbf{\nabla}\mathcal{L} (\mathbf{w} + \alpha\mathbf{p} )$ provides
$\mathbf{\nabla}\mathcal{L} (\mathbf{w}_{k+1})$ for the next iteration.
Note, however, that if the initial $\alpha_0$ in Alg.\;\ref{alg:interpLS}
satisfies (\ref{decrease}) and (\ref{curvature}), it is accepted
as the solution to (\ref{ls}) \emph{without} constructing and
minimizing $P_0(\alpha)$, regardless of the accuracy of
$\alpha_0$.

	\begin{algorithm}[tb!]
		\SetAlCapHSkip{0em}
		\small
		% 
% alg/interpLS.tex
% =================================================
% Author: Michael B Hynes, mbhynes@uwaterloo.ca
% License: GPL 3
% Creation Date: Wed 30 Sep 2015 04:26:48 PM EDT
% Last Modified: Fri 09 Oct 2015 01:07:22 AM EDT
% =================================================
\SetKwData{Left}{left}
\SetKwData{Up}{up}
\SetKwFunction{FindCompress}{FindCompress}
\LinesNumbered
\Indm
\KwIn{$\step_0 > 0$}
\KwOut{$\step_j$ satisfying \eqn{decrease} and
\eqn{curvature} for $\phi(\step)$}
\Indp
$I_0 \leftarrow [0,\step_0]$\\
$j \leftarrow 0$\\
\While{$\step_j$ does not satisfy \eqn{decrease} and \eqn{curvature}}{
	$[\step_l, \step_u] \leftarrow I_j$ \\
	$P_j(\alpha) \leftarrow$ interpolate
	$\phi(\step_l),\phi(\step_u),\dphi(\step_l),\dphi(\step_u)$\\
	$\step_{j+1} \leftarrow \arg\min_{\step} P_j(\step)$ \\
	$I_{j+1} \leftarrow$ update interval using $\step_l,\step_u,\step_{j+1}$  \\
	$j \leftarrow j + 1$
}
\KwRet $\step_j$

		\caption{\small Cubic Interpolating
	WA\ Line Search}
		\label{alg:interpLS}
	\end{algorithm}

\subsection{Polynomial Expansion Line Search.}
\label{sec:polyLS}
Our goal with the PELS\ method is to solve (\ref{ls})
accurately with as few \emph{distributed} operations as possible
when both $R (\mathbf{w})$ and 
$f(\mathbf{w};\mathbf{x}_i,y_i) $ are smooth with respect to $\mathbf{w}$. 
In this case, we consider the Taylor expansion of
$\phi(\alpha) = \mathcal{L} (\mathbf{w} + \alpha \mathbf{p} )$ in terms of $\alpha$,
which requires summing the 
Taylor expansions of each $f(\mathbf{w};\mathbf{x}_i,y_i) $ in (\ref{loss}) in addition to an
expansion for $R (\mathbf{w})$. % $(\x,\y)$ of the Taylor expansions of $\f$,
For an expansion about a step size $\alpha_j$, we have
\begin{equation*}\label{polyAll}
	\mathcal{L} (\mathbf{w} +\alpha\mathbf{p} ) = 
	Q_{\lambda} (\alpha)|_{\alpha_j} +
	\sum_{i=1}^n\sum_{\ell=0}^{\infty}
	b_{\ell,i}(\alpha-\alpha_j)^\ell
\end{equation*}
where $Q_{\lambda} (\alpha)|_{\alpha_j}$ is a polynomial
expansion of $\lambda R (\mathbf{w})$, and
${\left\{b_{\ell,i}\right\}}_{\ell = 0}^{\infty}$
are the coefficients in the expansion of $f(\mathbf{w};\mathbf{x}_i,y_i) $ that
depend explicitly on $(\mathbf{x}_i,y_i)$.
To make this representation amenable to a distributed
setting, we reorder the summations and
write the degree-$d$ approximation to $\mathcal{L} (\mathbf{w}+\alpha\mathbf{p} )$ as 
\begin{equation}\label{poly}
	W(\alpha; \mathbf{w},\mathbf{p} )|_{\alpha_j} =
		Q_{\lambda} (\alpha)|_{\alpha_j} + 
		\sum_{\ell=0}^{d}c_\ell(\alpha-\alpha_j)^\ell,
\end{equation}
which has a 
truncation error $\varepsilon_{d+1}(\alpha)|_{\alpha_j}$ as
\begin{equation}\label{polyApprox}
	W(\alpha;\mathbf{w},\mathbf{p} )|_{\alpha_j} = \mathcal{L} (\mathbf{w}+\alpha\mathbf{p} ) -
	\varepsilon_{d+1}(\alpha)|_{\alpha_j}.
\end{equation}
This error is an $\mathcal{O}\left[ (\alpha - \alpha_j)^{d+1} \right]$ term 
and hence small for $\alpha$ near $\alpha_j$.
Each coefficient $c_\ell$ in
(\ref{poly}) contains a summation over the dataset as
\begin{equation}\label{coeff}
	c_\ell = \sum_{i = 1}^{n} b_{\ell,i} = \sum_{i = 1}^{n}F_{\ell}(\mathbf{r},\mathbf{p} ;\mathbf{x}_i,y_i),
\end{equation}
where $\mathbf{r} = \mathbf{w} + \alpha_j\mathbf{p} $, and the functions
${\left\{F_\ell\right\}}_{\ell = 0}^{d}$ compute the coefficients
for the $\ell$th terms 
for a single observation.

The PELS\ algorithm minimizes the number of
distributed operations required to solve (\ref{ls})
by exploiting the following facts: (i) computing the coefficients
in (\ref{coeff})
is a \emph{parallelizable} operation, and (ii) once
the ${\left\{c_\ell\right\}}$ are computed,
$W(\alpha)|_{\alpha_j}$ is a useful approximation
in a neighbourbood of
$\alpha_j$.
The PELS\ method proceeds as follows.
Starting from the iterate $\alpha_j$, a polynomial
approximation $W(\alpha)|_{\alpha_j}$ is
constructed by computing the coefficients
${\left\{c_\ell\right\}}$ in
parallel, after which the subsequent iterate is determined
as
\begin{equation}\label{polyIterate}
	\alpha_{j+1} =
	\arg\min_{\alpha>0}W(\alpha;\mathbf{w},\mathbf{p} )|_{\alpha_j}.
\end{equation}
Solving (\ref{polyIterate}) is a subproblem that will generate further iterates in a
minimization routine, however, since the coefficients of
$W(\alpha)|_{\alpha_j}$ are \emph{fixed},
the minimization problem requires
no further distributed operations.
In addition, computing the first
and second derivatives $W^{\prime}|_{\alpha_j}$
and $W^{\prime\prime}|_{\alpha_j}$ are inexpensive
$\mathcal{O}\left[ d+1 \right]$ operations in a minimization routine
 when performed with
Horner's rule.
If the truncation error
$\varepsilon_{d+1}(\alpha_{j+1})|_{\alpha_j}$ in
(\ref{polyApprox}) is determined to be
too large, then $\alpha_{j+1}$ is likely a poor
approximation to $\alpha^*$, and the process repeats
iteratively: $\alpha_{j+1}$ is chosen as the new expansion
point, and the coefficients of
$W(\alpha)|_{\alpha_{j+1}}$ are computed.
This general form of the PELS\ algorithm is summarized in
Alg.\;\ref{alg:polyLS}, where the coefficients for $W(\alpha)|_{\alpha_j}$
are denoted by the vector 
$\mathbf{c} _j = [c_0,c_2,\ldots, c_{d}]^T$, and
the procedure {\sc calc\_coeffs} will
be described in \S\ref{sec:sparkLS}.
As a termination condition at step 7 of Alg.\;\ref{alg:polyLS},
we use an approximation to the fractional error in the value of
$W(\alpha_{j+1})|_{\alpha_j}$.
To estimate the term
$\varepsilon_{d+1}(\alpha_{j+1})|_{\alpha_j}$ in
the fractional error, 
we note that the truncation error in the
degree-$d$ Taylor polynomial is heuristically
bounded above by the error in the
polynomial of degree $d - 1$; that is, in a
relevant neighbourhood of the
expansion point $\alpha_j$,
$\left\lvert \varepsilon_{d}(\alpha) \right\rvert 
\approx \left\lvert c_{d}(\alpha - \alpha_j)^{d} \right\rvert
\geq \left\lvert \varepsilon_{d+1}(\alpha) \right\rvert$.
Hence, we approximate $\varepsilon_{d+1}(\alpha_{j+1})|_{\alpha_j}$ by
$c_{d}(\alpha_{j+1}
-\alpha_j)^{d}$ at step 5 of Alg.\;\ref{alg:polyLS}.

	\begin{algorithm}[tb!]
		\SetAlCapHSkip{0em}
		\small
		% 
% als/polyLS.tex
% =================================================
% Author: Michael B Hynes, mbhynes@uwaterloo.ca
% License: GPL 3
% Creation Date: Wed 30 Sep 2015 01:31:29 PM EDT
% Last Modified: Wed 20 Jan 2016 06:19:21 AM EST
% =================================================
\DontPrintSemicolon
\SetKwData{Left}{left}
\SetKwData{Up}{up}
\SetKwFunction{FindCompress}{FindCompress}
\Indm
\KwIn{$\w \in \Real^m,\p \in \Real^m,\step_0 > 0,
\tol > 0$}
\KwOut{$\step_j \approx \step^* = 
	\arg\min_{\step}\loss\left(
		\w + \step\p
	\right)$
}
\Indp
% \BlankLine
\nl $j \leftarrow 0$\\
% \nl $\errPoly(\step_j) \leftarrow \infty$\\
% \nl $\tol_j \leftarrow \infty$\\
% \nl \While{$\step_j$ is not sufficiently accurate}{
\nl \Repeat{\nl $\abs{\frac{\err_j}{\W(\step_j;\w,\p)\about{\step_{j-1}}}} \leq
	\tol$}{
% \nl	Compute $\Set{\polyC_i}$ for
% 	$\W(\step;\p,\w)|_{\step_j}$  \\
\nl	$\polyCoeff_j \leftarrow$ {\sc calc\_coeffs$(\w,\p,\step_j)$}\\% for $\W(\step;\p,\w)|_{\step_j}$\\
	% $\W(\step;\p,\w)|_{\step_j}$  \\
\nl	$\step_{j+1} \leftarrow \arg\min_{\step}\W(\step;\w,\p)|_{\step_j}$ \\
\nl	$\err_{j+1}\leftarrow$ approximate
$\errPoly(\step_{j+1})\about{\step_j}$ \\
\nl	$j \leftarrow j+1$ \\
}
\nl \KwRet $\step_j$
\BlankLine

\hrule
\BlankLine

\setcounter{AlgoLine}{0}
\SetKwProg{myproc}{Procedure}{}{}
\myproc{\sc calc\_coeffs$(\w,\p,\step_0)$}{
% \myproc{\sc calc\_coeffs$(\w,\p)$}{
\nl $\rr \leftarrow \w + \step_0 \p$\\
\nl Broadcast $\rr,\p$ to $\nnodes$ compute nodes\\
% \nl Broadcast $\w,\p$ to compute nodes\\
\nl\For{\normalfont compute node $t \in \Set{1,\ldots,\nnodes}$}{
\nl	\For{$\ell \in \Set{0,\ldots,\polyDegree}$}{
\nl		$\polyCoeffScalar^{[t]} \leftarrow \displaystyle 
			\sum_{\text{local }i_l} F_\ell(\rr,\p;\mathbf{x}_{i_l},y_{i_l})$
			% \sum_{\text{local }i_l} F_\ell(\w,\p;\mathbf{x}_{i_l},y_{i_l})$
	}
}
\nl $\polyCoeff_\lambda \leftarrow
	\frac{\lambda}{2}\left[\norm{\rr}^2_2,\;2\rr^T\p,\;\norm{\p}^2_2, 0, 0, \ldots\right]^T$ \\
\nl $\polyCoeff \leftarrow
	\polyCoeff_\lambda + \Ops_{t=1}^{\nnodes}\polyCoeff^{[t]} $\\
\nl \KwRet $\polyCoeff$\\
% $\step^* \leftarrow \arg\min_{\step}\W(\step;\wk,\pk)$
% \nl \KwRet $\Set{\polyCoeff}$
}

		\caption{\small Polynomial Expansion Line Search}
		\label{alg:polyLS}
	\end{algorithm}

\begin{figure}[bt!]
		\centering
		\resizebox{1.035\columnwidth}{!}{			\input{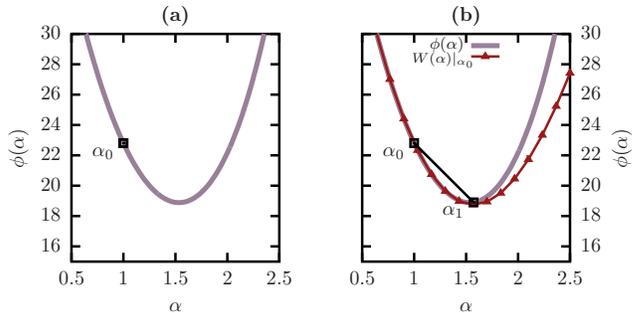}
		}
		\vspace{-0.9em}
		\caption{\emph{
	Example iterates ${\left\{\alpha_j\right\}}$ produced by
	\textbf{(a)} Alg.\;\ref{alg:interpLS} with $\nu_1 = 10^{-4}$ and
	$\nu_2
	= 0.9$ (implementation from \cite{more1994line}) and \textbf{(b)} 
	PELS\ with expansions using $d = 3$ for 
	$\phi(\alpha) = \alpha\,e^{\alpha} + e^{-(\alpha -4)}$.
}}
		\label{fig:ls_example}
	\end{figure}

When $\mathcal{L} (\mathbf{w})$ is itself polynomial, such as for
least squares regression and
low-rank matrix factorization as formulated by e.g.~Gemulla et
al.~\cite{gemulla2011large}, then $\varepsilon_{d+1}(\alpha) = 0$ for
sufficiently large $d$, and Alg.\;\ref{alg:polyLS}
can solve (\ref{ls}) \emph{exactly} in a single step without
further coefficient computations. 
For the more general case of analytic $\mathcal{L} (\mathbf{w})$, the
advantages of Alg.\;\ref{alg:polyLS} over Alg.\;\ref{alg:interpLS} are twofold.
Firstly, the communication costs are lesser
in the distributed operations: the
summation in (\ref{coeff}) to compute the
coefficients ${\left\{c_\ell\right\}}$
requires communicating $d+1$ \emph{scalar} values
for each $(\mathbf{x}_i,y_i)$, whereas
the distributed computation of $\mathbf{\nabla}\mathcal{L} $
in Alg.\;\ref{alg:interpLS} communicates the $m$-dimensional gradient
vectors ${\left\{\mathbf{\nabla} f(\mathbf{w};\mathbf{x}_i,y_i)\right\}}$.
The second---and most important---reason is that, unlike
standard line searches that compute values for $\mathcal{L} $
and $\mathbf{\nabla}\mathcal{L} $ for only a \emph{single}
$\alpha_j$, the distributed
operation computing the ${\left\{c_\ell\right\}}$ in PELS\
produces a model which is valid for 
an entire neighbourhood of $\alpha_j$;
thus, if $\alpha_j$ is close 
to $\alpha^*$, the PELS\ method can compute a very accurate
approximation to $\alpha^*$ with only a single
$\mathcal{O}\left[ n \right]$ pass over the dataset to evaluate the polynomial
coefficients.
This is illustrated in Fig.\;\ref{fig:ls_example},
where both Alg.\;\ref{alg:interpLS} and Alg.\;\ref{alg:polyLS} have been applied
to the function  $\phi(\alpha) = \alpha\,e^{\alpha} +
e^{-(\alpha -4)}$. 
In Fig.\;\ref{fig:ls_example}a, iterates produced by Alg.\;\ref{alg:interpLS} (implemented as
in \cite{more1994line} with $\nu_1 = 10^{-4}$ and $\nu_2 = 0.9$)
are shown, however only a single iterate has been generated since
(\ref{decrease}) and (\ref{curvature}) are
satisfied at the input $\alpha_0 = 1$; as
such, the algorithm terminates with a relatively inaccurate
minimum, rather than
compute $\mathbf{\nabla}\mathcal{L} $ with an $\mathcal{O}\left[ n \right]$ pass over the data
for another step $\alpha_1$.
On the other hand, Fig.\;\ref{fig:ls_example}b shows PELS\  
with a degree-3 polynomial approximation $W|_{\alpha_0}$ to
$\phi(\alpha)$ for the same $\alpha_0 = 1$.
Here, $\alpha_1$ is determined as the solution to
(\ref{polyIterate}), and is visibly more accurate than
$\alpha_0$ in Fig.\;\ref{fig:ls_example}a (only $d = 3$ is
shown since $W|_{\alpha_0}$ could not
be distinguished from $\phi(\alpha)$ at this scale for larger degrees).
This example is not a toy problem, but a case that we have often observed experimentally for the LBFGS\
algorithm: though the initial step size $\alpha_0 = 1$ is
frequently accepted in Alg.\;\ref{alg:interpLS}, $\alpha_0$ is an
inaccurate
solution to (\ref{ls}).

In our implementation of PELS, we
used the NR\ method as the optimization routine for
minimizing the polynomial $W|_{\alpha_j}$ in
step 4 of Alg.\;\ref{alg:polyLS},
which is equivalent to solving for the roots of
$W^{\prime}|_{\alpha_j}$.
However, it is possible
that $W^{\prime}|_{\alpha_j} < 0$ for $\alpha>0$ if $\alpha_j$ is far from $\alpha^*$ (i.e.~$W^{\prime}|_{\alpha_j}$ has no relevant real roots).
In this case, instead of using a minimization routine,
$\alpha_{j+1}$ was determined by 
NR\ iteration at $\alpha_j$ as
\begin{equation}\label{newton}
		\alpha^{NR}
		= \alpha_j-\frac{\phi^{\prime}(\alpha_j)}{\phi^{\prime\prime}(\alpha_j)} 
		= \alpha_j-\frac{c_1}{2c_2},
\end{equation}
where only the coefficients $c_1$ and $c_2$
appear since $\phi(\alpha_j) = W(\alpha_j)|_{\alpha_j}$.
Since in practice, NR iteration converged within machine
precision to a stationary point
of $W|_{\alpha_j}$ in very few iterations,
(\ref{newton}) was used as a default whenever
the NR method starting from $\alpha_j$ failed to
compute a zero gradient within a tolerance of $10^{-15}$
in fewer than 10 iterations.
Equation (\ref{newton}) also indicates that very fast convergence
is expected for PELS\ when $\alpha_j$ is close to 
$\alpha^*$, since
NR iteration is equivalent to using only a $d=2$
approximation, whereas the PELS\ method can utilize
arbitrary $d > 2$.
Finally, note that NR\ iteration for step 4 may not be convergent for some initial
guesses. Naturally, 
other minimization algorithms (e.g.~bisection or
backtracking) may be applied to
$W|_{\alpha_j}$, however empirically we have
found that a simple step size selection method (described in
\S\ref{sec:implementation}) to generate 
guesses for  $\alpha_0$ in each invocation of
PELS\ produces stable results.

\subsection{Logistic Regression.}
Consider a logistic regression model with
$L_2$-regularization for 
binary classification of each $\mathbf{x}_i \in \mathbbm{R}^m$ with a group label $y_i \in
{\left\{0,1\right\}}$.  The loss function is derived from the
maximum log likelihood of the probability 
$\text{Pr} \left[ y_i = 0|\mathbf{x}_i,\mathbf{w} \right]$ as
\begin{equation}\label{logisticLoss}
	\mathcal{L} (\mathbf{w}) = \frac{\lambda}{2}\|\mathbf{w}\|_2^2 + \frac{1}{n}\sum_{i=1}^n \left(
		\log[e^{-\mathbf{w}^T\mathbf{x}_i} + 1] + \mathbbm{N}(y_i) \mathbf{w}^T\mathbf{x}_i
	\right)
\end{equation}
where $\mathbbm{N}(y_i) = 1 - \mathbbm{I}(y_i)$ and $\mathbbm{I}(y_i)$ is a
boolean indicator function equal to 1 only if $y_i \neq 0$.
It can be shown that (\ref{logisticLoss}) is strictly convex
and has a unique global minimizer~\cite{lin2008trust}.
Since (\ref{logisticLoss}) is infinitely differentiable, we
may use a Taylor expansion about $\alpha_0$ for
$\phi(\alpha) = \mathcal{L} (\mathbf{w} + \alpha \mathbf{p} )$.
Denoting the ray $\mathbf{w} + \alpha_0 \mathbf{p} $ by $\mathbf{r}$, and using the 
simplifications $p_i = \mathbf{p} ^T\mathbf{x}_i$ and $r_i = e^{-\mathbf{r}^T\mathbf{x}_i}$,
the expansion up to 4th order is
\begin{align}\label{logisticPoly}
		\phi(\alpha)
		&= \frac{1}{n}\sum_{i=1}^n \left( \log
			\left[r_i+1\right] + \mathbbm{N}(y_i)\,\mathbf{r}^T\mathbf{x}_i
		\right) + \frac{\lambda}{2}\|\mathbf{r}\|_2^2 \\ 
		&+ (\alpha-\alpha_0) \left[
			\frac{1}{n}\sum_{i=1}^n p_i \left(
			\mathbbm{N}(y_i)-\frac{r_i}{r_i+1}\right)
			+ \lambda	\mathbf{p} ^T\mathbf{r} \right] \nonumber\\ 
		&+ (\alpha-\alpha_0)^2 \left[
			\frac{1}{n}\sum_{i=1}^n \frac{p_i^2r_i}{2(r_i+1)^2}
			+ \frac{\lambda}{2}\|\mathbf{\mathbf{p} }\|_2^2 
			\right] \nonumber\\
			&+ \frac{(\alpha-\alpha_0)}{n}^3\sum_{i=1}^n
			{{p_i^3r_i\left(r_i-1\right)}\over{6\left(r_i+1\right)^3}} \nonumber\\
		& + \mathcal{O}\left[ (\alpha-\alpha_0)^4 \right].  \nonumber
\end{align}
The terms $p_i$ and $r_i$ in coefficients of
(\ref{logisticPoly}) can be computed in parallel for each
instance, $(\mathbf{x}_i,y_i)$.
Furthermore, higher order coefficients in (\ref{logisticPoly}) depend
only on the scalars $p_i$ and $r_i$.
Therefore, the only vector
operations required in parallel are $\mathbf{p} ^T\mathbf{x}_i$ and $\mathbf{r}^T\mathbf{x}_i$.

\subsection{Apache Spark.}
\label{sec:spark}

Apache Spark is a fault-tolerant, in-memory cluster
computing framework designed to supersede MapReduce by
maintaining program data in memory as much as possible
between distributed operations.
The Spark environment is built upon two components: a data
abstraction, termed a resilient distributed dataset (RDD)
\cite{zaharia2012resilient}, and the task scheduler, which
uses a \emph{delay scheduling} algorithm
\cite{zaharia2010delay}.
A Spark cluster is composed of a set of (slave) executor
programs running on $n_p$ compute
nodes, and a master program running on a master node
that is
responsible for scheduling and allocating tasks 
to the compute nodes based on data locality.
We describe the fundamental aspects of RDDs and the scheduler below.

RDDs are immutable, distributed datasets that are evaluated
lazily via their provenance information---that is, their
functional relation to other RDDs or datasets in stable
storage.  To describe an RDD, consider an immutable
distributed dataset $D$ of $n$ records with homogeneous
type: $D = \bigcup_i^n d_i$ with $d_i \in \mathcal{D}$.  The
distribution of $D$ across a computer network of $n_p$
nodes ${\left\{v_{t}\right\}}$, such that $d_i$ is stored in
memory or on disk on node $v_{t}$, is termed its
\emph{partitioning} according to a partition function
$P(d_i) = v_{t}$.  If $D$ is expressible as a finite
sequence of deterministic operations on other datasets
$D_1,\ldots, D_l$ that are either RDDs or persistent
records, then its lineage may be written as a directed
acyclic graph $L$ formed with the parent datasets
${\left\{D_l\right\}}$ as the vertices, and the operations along the
edges.  Thus, an RDD of type $\mathcal{D}$ (written
$\mathtt{\small{RDD}}\left[  \mathcal{D} \right] $) is the tuple $(D,P,L)$.

Physically computing the records ${\left\{d_i\right\}}$ of an RDD is
termed its \emph{materialization}, and is managed by the
Spark scheduler program.  To allocate computational tasks to
the compute nodes, the scheduler traverses an RDD's lineage
graph $L$ and divides the required operations into
\emph{stages} of local computations on parent RDD partitions.
Suppose that $R_0 = (\bigcup_i x_i,P_0,L_0)$ were an
RDD of numeric type $\mathtt{\small{RDD}}\left[  \mathbbm{R} \right] $, and let $R_1 =
(\bigcup_i y_i, P_1, L_1)$ be the RDD resulting from
the application of a function
$f:\mathbbm{R}\rightarrow\mathbbm{R}$ to each record of $R_0$.
To compute ${\left\{y_i\right\}}$, $R_1$ has only a single parent in
the graph $L_1$, and hence the set of tasks to
perform is ${\left\{f(x_i)\right\}}$.  This type of operation is termed
a \emph{map} operation, and has a \emph{narrow} lineage
dependency: $P_1 = P_0$, and  the scheduler would
allocate the task $f(x_i)$ to a node that stores $x_i$ since each $y_i$
may be computed locally from $x_i$.

Stages consist only of local map operations, and are bounded
by \emph{shuffle} operations that require communication and
data transfer between the compute nodes.  For example,
shuffling is necessary to perform \emph{reduce} operations
on RDDs, wherein a scalar value is produced from an
associative binary operator applied to each element of the
dataset.  In implementation, a shuffle is conducted by
writing the results of the tasks in the preceding stage,
${\left\{f(x_i)\right\}}$, to a local file buffer. 
These shuffle files may or may not be written to disk,
depending on the operating system's page table, and are
fetched by remote nodes as needed in the subsequent stage
via a TCP connection. %~\cite{ousterhout2015shuffle}.
Shuffle file fetches occur asynchronously, and multiple
connections between compute nodes to transfer information
are made in parallel.  
In addition, map tasks on the previous stage's
results that are stored locally by a
compute node will occur concurrently with remote fetches.

\begin{figure}[bt!]
		\centering
		\resizebox{0.71\columnwidth}{!}{			\input{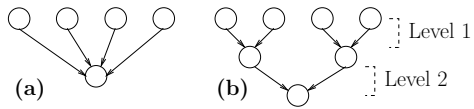}
		}
		\vspace{-0.9em}
		\caption{\emph{Comparison of an (a) all-to-one
	communication where every node, depicted by circles,
	communicates results with the
	driver (bottom circle), and a (b) multi-level scheme with
	$n_l = 2$, such that the final reduction to the
	driver requires communication with only 2 nodes.
}}
		\label{fig:aggregate}
	\end{figure}

A reduce operation on an RDD of type $\mathtt{\small{RDD}}\left[  \mathcal{A} \right] $ produces a
scalar value of type $\mathcal{A}$ by an application of an associative
binary operator $\hspace{1pt}\mathbf{\oplus}\hspace{1pt}:\mathcal{A}\times\mathcal{A}\rightarrow\mathcal{A}$ to all of the
elements ${\left\{a_i\right\}}$ as
$\bigoplus_{i=1}^n a_i = a_1 \hspace{1pt}\mathbf{\oplus}\hspace{1pt} a_2 \hspace{1pt}\mathbf{\oplus}\hspace{1pt} \cdots \hspace{1pt}\mathbf{\oplus}\hspace{1pt} a_n$.
Reduce operations require first performing \emph{local}
reductions on the partitions stored locally by each node before
communicating the nodes' results to the driver.
The local reduction of the partition of ${\left\{a_i\right\}}$ stored on node
$v_t$ is written $a^{[t]}$ and computed as
$a^{[t]} = \bigoplus_{i_l \in \mathcal{I}_t} a_{i_l}$ for $\mathcal{I}_t =
{\left\{i:P(a_{i}) = v_t\right\}}$
and requires no communication between nodes.
The full reduction for multiple nodes is then
$\bigoplus_{t=1}^{n_p}
a^{[t]}$, which incurs communication costs
dependent on $n_p$ and the size (in bytes) of
an element in $\mathcal{A}$: reducing
scalars ($\mathcal{A} = \mathbbm{R}$) is cheaper than reducing vectors
($\mathcal{A} = \mathbbm{R}^p$).
Additionally, reduce operations on an RDD may be performed in two ways: either through an all-to-one communication pattern
(Fig.\;\ref{fig:aggregate}a) in which all $n_p$ local results
from the compute nodes are communicated to the host machine
on which the driver program is running, or through a
multi-level tree communication \cite{agarwal2014reliable}
where intermediary (locally reduced) results are aggregated by
compute nodes in $n_l$ levels before being
transferred
to the driver \cite{yavuz2014mllib} (Fig.\;\ref{fig:aggregate}b).

\section{PELS\ Implementation \& Performance Tests.}
\label{sec:methods}

\subsection{PELS\ Implementation In Apache Spark.}
\label{sec:sparkLS}

The PELS\ method in Alg.\;\ref{alg:polyLS} was implemented in
Apache Spark based on the LBFGS\
implementation in the Spark 1.5 codebase.
The observations $\mathcal{R} =
{\left\{(\mathbf{x}_i,y_i)\right\}}$ were stored in an RDD with type 
$\mathtt{\small{RDD}}\left[  (\mathbbm{R}^p,\mathbbm{R}) \right] $, where the vectors ${\left\{\mathbf{x}_i\right\}}$
were either sparse or dense vectors with double precision.
The parameter and search direction vectors ${\left\{\mathbf{w}_{k}\right\}}$ and
${\left\{\mathbf{p} _k \right\}}$ were not partitioned over the compute nodes,
but stored as dense vector
objects on the master node.
Communicating the vectors $\mathbf{p} _k $ and $\mathbf{r}_k = \mathbf{w}_{k} + \alpha_j\mathbf{p} _k $ to the compute nodes in the
cluster for any $\alpha_j$ in the line search
was performed via a torrent broadcast~\cite{yavuz2014mllib}.
To compute the coefficient vector $\mathbf{c} _j$ in any PELS\
iteration, the functions
${\left\{F_{\ell}(\mathbf{r}_k,\mathbf{p} _k ;\mathbf{x}_i,y_i)\right\}}_{\ell=0}^{d}$ in
(\ref{coeff}) were
calculated as parallel map operations on the RDD of $\mathcal{R}$, and the
resulting $\mathbf{c} _j \in
\mathbbm{R}^{d+1}$ was summed via a multi-level reduce
operation (as in Fig.\;\ref{fig:aggregate}b), where the operator $\hspace{1pt}\mathbf{\oplus}\hspace{1pt}$ was
vector addition.
The $L_2$ regularization terms
were computed by the master node 
as the vector
$\mathbf{c} _{\lambda} = \frac{\lambda}{2}[\|\mathbf{r}_k\|_2^2,
2\mathbf{r}_k^T\mathbf{p} _k , \|\mathbf{p} _k \|_2^2, 0,0,\ldots]$ (with zeros appended
so $\mathbf{c} _{\lambda}
\in \mathbbm{R}^{d+1}$) and added to $\mathbf{c} _j$.
This is detailed in
procedure {\sc compute\_coeffs}
in Alg.\;\ref{alg:polyLS}, in which the local and global reductions
occur at steps 5 and 7, respectively.

\subsection{PELS\ for LBFGS, NCG, \& GD.}
\label{sec:implementation}
The PELS\ algorithm with $d = 5$
was used to accelerate the following algorithms in Apache Spark
for unconstrained optimization of smooth
$L_2$-regularized loss
functions: %with $\ell_2$ regularization:
(i) LBFGS,
(ii) NCG,
and (iii) GD.
The implementations using PELS\ will be henceforth
denoted with a suffix \emph{-P} as LBFGS-P, NCG-P, and
GD-P.
For LBFGS-P, NCG-P, and GD-P, once the search direction $\mathbf{p} _k $ was determined at $\mathbf{w}_{k}$ in
the $k$th iteration, Alg.\;\ref{alg:polyLS} was used to solve (\ref{ls}). 
For both the WA\ line search in LBFGS\ and PELS\ in LBFGS-P, the initial step size  
in each invocation
(i.e.~the
input $\alpha_0$ in Alg.\;\ref{alg:interpLS} and Alg.\;\ref{alg:polyLS})
was taken to be $1$, which is the
recommended trial step for the
LBFGS\ algorithm \cite{liu1989limited}.
The NCG and GD algorithms with and without PELS\
both used a 
standard scaling formula \cite{nocedal2006numerical}
to determine the initial step size $\alpha_k^{0}$ for the
line search invocation in the $k$th
iteration as
\begin{equation*}\label{initStep}
	\alpha_k^0 = \alpha_{k-1}
	\frac{\mathbf{\nabla}\mathcal{L} (\mathbf{w}_{k-1})^T\mathbf{p} _{k-1}}{\mathbf{\nabla}\mathcal{L} (\mathbf{w}_{k})^T\mathbf{p} _k },
\end{equation*}
where $\alpha_{k-1}$ is the
accepted step size from the previous iteration as in
(\ref{update}), and $\alpha^0_{k} = 1$ for $k=0$.
For NCG-P and GD-P, the search directions used
to compute the coefficients in PELS\ and
subsequently update $\mathbf{w}_{k}$
 were normalized as $\mathbf{p} _k  \leftarrow \mathbf{p} _k  /
\|\mathbf{p} _k \|_2$.
The NCG-P update rule was the positive Polak-Ribi\`ere
formula \cite{gilbert1992global}, with a  
Powell restart condition \cite{powell1977restart}.
For LBFGS-P, the initial inverse Hessian approximation in each
iteration used Liu and Nocedal's M3 scaling
\cite{liu1989limited},
as in the LBFGS\ code for Spark version 1.5.

\subsection{Algorithm Performance Tests.}
\label{sec:tests}
To evaluate the efficacy of PELS, performance tests were conducted with a
logistic regression model as in
(\ref{logisticLoss}) that compared
the LBFGS-P, NCG-P, and GD-P implementations with
LBFGS, NCG, and GD using
the WA\ 
cubic interpolating line search \cite{al1986efficient,nocedal2006numerical}
in the Breeze library~\cite{breeze} with $\nu_1 =
10^{-4}$ and $\nu_2 = 0.9$, which is used in Spark's LBFGS\
implementation.
The LBFGS, NCG, and GD algorithms using a WA\
line search were implemented such that all vector operations
other than the line search were identical to those in the LBFGS-P, NCG-P and
GD-P implementations, respectively\footnote{
	We found that Spark's LBFGS\ code had
	significant overheads stemming from the 
	deeply nested data structures
	in Breeze; we wrote the LBFGS\
	algorithm in an imperative style that took only $\sim$65\%
	as much time as the Spark code on large problems.
}.
The Powell restart threshold in NCG-P was 0.2, however NCG
had a restart threshold of 1.0 since
smaller values often triggered restarts in \emph{every} iteration,
reverting the method to GD.
For both
LBFGS\ and LBFGS-P,
a history of 5 corrections was used in our tests,
as recommended for large problems
\cite{liu1989limited}.

The logistic regression models were trained on 
binary classification datasets procured from
the LIBSVM repository
\cite{chang2011libsvm}, which are listed in
in Tab.\;\ref{tab:datasets}. 
This table also gives summary information about the
datasets' respective sizes of $n$ and $m$,
mean number of nonzero elements in $\mathbf{x}_i$, 
ratio of the number of true (nonzero) $y_i$ labels
to false $y_i$ labels as $n_+:n_-$, 
and magnitudes of $\lambda$ used in
(\ref{logisticLoss}).
The dense ${\left\{\mathbf{x}_i\right\}}$ in the {\sc e}psilon\, dataset were represented
by contiguous dense vectors, while the other datasets'
instances were stored in compressed sparse vector format.
All $\mathbf{x}_i$ were further augmented as $\mathbf{x}_i^T \leftarrow[\mathbf{x}_i^T\;
1]^T$ to implicitly include a
constant offset term in the inner product $\mathbf{w}^T\mathbf{x}_i$.
For all datasets, each algorithm 
was run with $\mathbf{w}_0 = \mathbf{0}$.
In each iteration, 
computed values of $\mathcal{L} (\mathbf{w}_{k})$ and
$\|\mathbf{\nabla}\mathcal{L} (\mathbf{w}_{k})\|_2$ were written to
the standard output filestream. 
The values of $\theta$ used in Alg.\;\ref{alg:polyLS} were $\theta =
10^{-4}$ (approximately 0.01\% error), except for the
{\sc e}psilon\ and {\sc kdd-a}\ datasets, on
which $10^{-6}$ was used. While smaller $\theta$ generated more
accurate step sizes, the additional PELS\
iterations increased the computational time; values of
$\theta \in [10^{-6},10^{-4}]$ were a good compromise
between accuracy and speed.

	\begin{table}[t!]
		\centering
		\caption{\emph{
	Properties of LIBSVM classification datasets used in
	numerical experiments, and respective magnitudes of
	$\lambda$. 
}}
		\label{tab:datasets}
		\footnotesize
		\begin{tabular}{lrrccc}\hline\hline
Dataset&$n$&$m$&$\nnz{\x}$&$\ratio{n_+}{n_-}$&$\lambda$\\
\hline
\eps &400,000 &2,001 &2001 &\ratio{1.0}{1.0} &$10^{-8}$ \\
{\sc rcv\scriptsize{1}} (test) &677,399 &47,237 &$ 74 \pm 54$ &\ratio{1.1}{1.0} &$10^{-7}$ \\
\urldata &2,396,130 &3,231,962 &$ 117 \pm 17$ &\ratio{1.0}{2.0} &$10^{-8}$ \\
\kdda &8,407,752 &20,216,831 &$ 37 \pm 9$ &\ratio{5.8}{1.0} &$10^{-9}$ \\
\kddb &19,264,097 &29,890,095 &$ 29 \pm 8$ &\ratio{6.2}{1.0} &$10^{-9}$ \\
\hline\hline\end{tabular}

	\end{table}

All performance tests were performed on a computing cluster
composed of: 16
homogeneous compute nodes, 1 storage node hosting a network
filesystem, and 1 master node.
The nodes were interconnected by a 10 Gb ethernet
managed switch (PowerConnect 8164; Dell).
Each compute node was a 64 bit rack server (PowerEdge R620;
Dell) running Linux kernel 3.13 with two 8-core 2.60 GHz
processors (Xeon E5-2670; Intel) and 200 GB of SDRAM.
The master node had identical processor specifications and
512 GB of RAM.
Compute nodes were equipped with six ext4-formatted 600 GB
SCSI hard disks, each with 10,000 RPM nominal speed.
The storage node (PowerEdge R720; Dell)
contained two 6-core 2 GHz processors (Xeon E5-2620; Intel),
64 GB of memory, and a hard drive speed of
7,200 RPM.

Our Apache Spark assembly was built from a snapshot of the
version 1.5 master branch using Oracle's Java 7
distribution, Scala 2.10, and Hadoop 2.0.2.  Input files to
Spark programs were stored on the storage node in plain
text.
The compute nodes' local SCSI drives were used for both Spark spilling directories 
and Java temporary files.
Shuffle files were consolidated into larger files, as
recommended for ext4 filesystems \cite{davidson2013optimizing},
and Kryo serialization was used.
In our experiments, the Spark driver was executed on the
master node in standalone client mode, and a single instance
of a Spark executor was created on each compute node.
All RDDs had 256 partitions,
corresponding to 1 partition per available physical core;
performance decreased in general as more partitions were
used.
Finally, $n_l$ was set to
$\log_2{16}$ in all tree aggregations; empirically, this was faster for large datasets than the
default of $n_l = 2$.

\section{Results \& Discussion.}
\label{sec:results}
Our performance results are presented in two parts. 
In our initial tests, we are interested purely in the
convergence improvements possible with PELS, and hence
consider the {\sc e}psilon\ and {\sc rcv\footnotesize{1}}\ datasets
that have little noise and are well-posed with $m \ll n$.
In these tests, we present high-accuracy convergence
traces since
$\mathcal{L} (\mathbf{w}^*)$ may be computed to machine precision.
By contrast, real-life machine learning
applications are often ill-posed and have statistical
errors in
the model or instances that exceed the numerical optimization
error, such that $\mathbf{w}^*$ is computed to only a few significant
digits until the training loss ceases to
decrease appreciably \cite{bousquet2008tradeoffs}.
To demonstrate that PELS\ is effective in this practical
setting as well, we present results on the large and
ill-posed {\sc url}, {\sc kdd-a}, and {\sc kdd-b}\ datasets, and compute
the acceleration in reaching terminal values of the
training label prediction accuracy, $\exp\{-\mathcal{L} (\mathbf{w}_{k})\}$, achievable by using PELS.

\begin{figure}[t!]
		\centering
		\resizebox{1\columnwidth}{!}{			\input{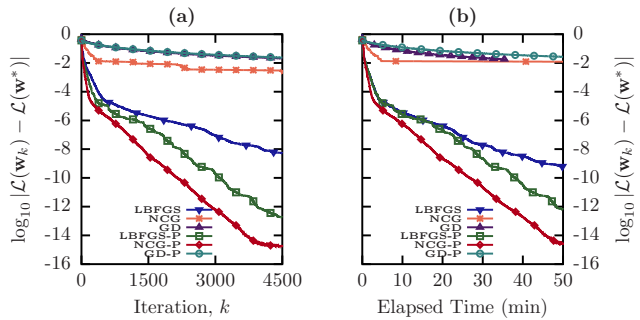}
		}
		\vspace{-0.9em}
		\caption{\emph{
	Convergence traces in regularized loss for the {\sc e}psilon\
	dataset in \textbf{(a)} iterations and \textbf{(b)}
	elapsed time.
}}
		\label{fig:eps-trace}
	\end{figure}

For the {\sc e}psilon\ and {\sc rcv\footnotesize{1}}\ datasets, respectively,
Fig.\;\ref{fig:eps-trace} and Fig.\;\ref{fig:rcv1-trace} show the traces of 
$\left\lvert \mathcal{L} (\mathbf{w}_{k})-\mathcal{L} (\mathbf{w}^*) \right\rvert$ as a function of {\bf (a)}\
iterations and {\bf (b)}\ clock time for the algorithms
considered in \S\ref{sec:tests}.
In both plots, $\mathcal{L} (\mathbf{w}^*)$ was determined as the
minimal loss computed by any algorithm within the maximum number of
iterations.
That $\mathbf{w}^*$ was computed accurately is evinced
by the gradient norm at $\mathbf{w}^*$: $\|\mathbf{\nabla}\mathcal{L} (\mathbf{w}^*)\|_2$
was $1.3 \times 10^{-11}$ for the {\sc e}psilon\ dataset
and $8.7 \times 10^{-12}$ for the {\sc rcv\footnotesize{1}}\ dataset, as computed by NCG-P.
It is notable in both plots that NCG-P has drastically
outperformed both the NCG and LBFGS\ algorithms that use a standard 
WA\ line search.
LBFGS-P outperforms LBFGS\ in iterations and time for
the {\sc e}psilon\ dataset, and performs similarly to LBFGS\
on the {\sc rcv\footnotesize{1}}\ dataset.
GD-P and GD are virtually indistinguishable at
the scale of Fig.\;\ref{fig:eps-trace} and Fig.\;\ref{fig:rcv1-trace} since
there is little substantive difference in the two
algorithms' traces.

	\begin{figure}[t!]
		\centering
		\resizebox{1\columnwidth}{!}{			\input{rcv1-trace}
		}
		\vspace{-0.9em}
		\caption{\emph{
	Convergence traces in regularized loss for the {\sc rcv\footnotesize{1}}\
	dataset in \textbf{(a)} iterations and \textbf{(b)}
	elapsed time.
}}
		\label{fig:rcv1-trace}
	\end{figure}

	\begin{figure}[t!]
		\centering
		\resizebox{1\columnwidth}{!}{			\input{url-trace}
		}
		\vspace{-0.9em}
		\caption{\emph{
	Convergence traces in regularized loss for the {\sc url}\
	dataset in \textbf{(a)} iterations and \textbf{(b)}
	elapsed time.
}}
		\label{fig:url-trace}
	\end{figure}

	\begin{figure}[t!]
		\centering
		\resizebox{1\columnwidth}{!}{			\input{kdda-trace}
		}
		\vspace{-0.9em}
		\caption{\emph{
	Convergence traces in regularized loss for the {\sc kdd-a}\
	dataset in \textbf{(a)} iterations and \textbf{(b)}
	elapsed time.
}}
		\label{fig:kdda-trace}
	\end{figure}

	\begin{figure}[t!]
		\centering
		\resizebox{1\columnwidth}{!}{			\input{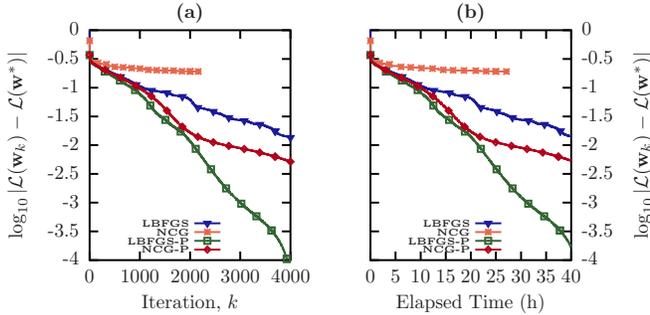}
		}
		\vspace{-0.9em}
		\caption{\emph{
	Convergence traces in regularized loss for the {\sc kdd-b}\
	dataset in \textbf{(a)} iterations and \textbf{(b)}
	elapsed time.
}}
		\label{fig:kddb-trace}
	\end{figure}

\begin{figure}[t!]
		\centering
		\resizebox{1\columnwidth}{!}{			\input{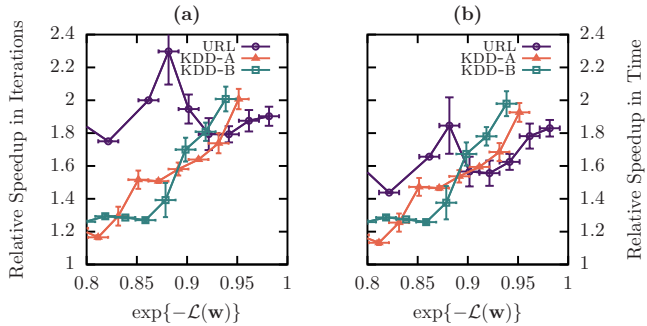}
		}
		\vspace{-0.9em}
		\caption{\emph{
	Speedup of LBFGS-P over LBFGS\ as a function of
	approximate training accuracy, computed for
	{\bf (a)}\ iterations and {\bf (b)}\ time.
}}
		\label{fig:bfgs_speedup}
	\end{figure}

The convergence traces for the large {\sc url}, {\sc kdd-a}, and {\sc kdd-b}\
datasets are shown in
Fig.\;\ref{fig:url-trace}, 
Fig.\;\ref{fig:kdda-trace},
and Fig.\;\ref{fig:kddb-trace}, respectively (traces for GD and GD-P have been omitted from these
plots since these algorithms made little progress towards
the solution).
These datasets required considerably more iterations and
training time, and the LBFGS-P and NCG-P algorithms have
outperformed
their counterparts by multiple decimal digits in accuracy. 
However, since in many machine learning applications,
the training procedure is halted once 
$\exp\{-\mathcal{L} (\mathbf{w}_{k})\}$ reaches a
plateau, only $\log_{10}\left\lvert \mathcal{L} (\mathbf{w}_{k}) - \mathcal{L} (\mathbf{w}^*) \right\rvert
\approx -3$ may be necessary in practice.
Bearing this,
the speedup factors as a function of training
accuracy for LBFGS-P relative to LBFGS\ were computed for the {\sc url}\
and both {\sc kdd} datasets, and are shown in
Fig.\;\ref{fig:bfgs_speedup}, where the factors in
Fig.\;\ref{fig:bfgs_speedup}a have been computed using the
iterations required and the factors in
Fig.\;\ref{fig:bfgs_speedup}b have been computed for clock time.
These speedup factors were determined by finding the first iterate
produced by LBFGS-P that reached the same or greater value
of $\exp\{-\mathcal{L} (\mathbf{w}_{k})\}$ as LBFGS,
and the error bars in this plot show the standard deviation
about the mean speedup ratios computed in non-overlapping windows of width 0.02
along the x-axis.
To reach the terminal training accuracies (e.g.\;$\sim$97\% for 
{\sc url}),
Fig.\;\ref{fig:bfgs_speedup} shows factors of 1.8--2 in both iterations
and clock time are achievable on the {\sc kdd-a}, {\sc kdd-b}, and
{\sc url}\ problems.

To complement the convergence traces and explain the speedup
in clock time
of the PELS\ method in a distributed setting,
Tab.\;\ref{tab:all} presents timing measurements for the distributed
operations, averaged for each algorithm over all
iterations.
The quantity $\tau_{\text{fg}}$ represents the mean
clock time required to evaluate $\mathcal{L} $ and $\mathbf{\nabla}\mathcal{L} $
through a map operation on the RDD of $\mathcal{R}$ with a
subsequent aggregation.
However, since $\mathcal{L} $ does not need to be computed
explicitly in the PELS\ method, $\tau_{\text{fg}}$ denotes the time
required to compute \emph{only} $\mathbf{\nabla}\mathcal{L} $ for LBFGS-P, NCG-P, and
GD-P (although evaluating $\mathbf{\nabla}\mathcal{L} $ alone takes as much
time as evaluating $\mathcal{L} $ and $\mathbf{\nabla}\mathcal{L} $
simultaneously). 
For the algorithms using the PELS\ method, $\tau_{\ell\text{s}}$ gives the
mean time to both compute and aggregate the coefficient
vector
$\mathbf{c} _j$.
The quantity $n_e$ represents the average number of
iterations per line search in the respective manners in
which they were conducted: for the
WA\ algorithms, $n_e$ is the mean number
of function/gradient evaluations per line search, and
for the PELS\ algorithms, it represents the mean number of
coefficient evaluations performed in Alg.\;\ref{alg:polyLS} per line search.
Thus, in each outer iteration of the optimization method,
PELS\ algorithms perform 1 gradient
computation taking $\tau_{\text{fg}}$ seconds and then $n_e$
operations lasting $\tau_{\ell\text{s}}$ seconds, while
WA\ algorithms perform $n_e$
operations taking $\tau_{\text{fg}}$ seconds.
The total time per iteration is shown as $\tau_{\text{tot}}$.
All uncertainty bounds show the standard deviation about the
mean value; 
no uncertainty is given for $n_e$ since it was taken as the ratio of the total
number of line search calls to outer iterations.

	\begin{table}[t!]
		\centering
		\caption{\emph{
	Mean clock times for each iteration ($\tau_{\text{tot}}$), evaluating $\mathcal{L} $ and
	$\mathbf{\nabla}\mathcal{L} $ ($\tau_{\text{fg}}$), and computing $\mathbf{c} _j$
	($\tau_{\ell\text{s}}$). The average number of function calls or line
	search iterations per outer iteration of the optimization
	algorithms are given by $n_e$.
}}
		\label{tab:all}
		\footnotesize
		\scriptsize
\begin{tabular}{llcccr}
\hline
\hline
 &Alg. &$\iterTime$ (s) &$\costTime$ (s) &$\polyEvalTime$ (s) &$\numEvals$ \\
\hline
\eps& & & & \\
\hline
& LBFGS&$0.5\pm0.2$&$0.42\pm0.07$&&1.15\\
& NCG&$4\pm2$&$0.43\pm0.08$&&9.92\\
& GD&$0.42\pm0.10$&$0.42\pm0.08$&&1.01\\
& LBFGS-P&$0.7\pm0.1$&$0.42\pm0.07$&$0.32\pm0.06$&1.00\\
& NCG-P&$0.7\pm0.1$&$0.42\pm0.07$&$0.32\pm0.06$&1.00\\
& GD-P&$0.7\pm0.1$&$0.42\pm0.07$&$0.32\pm0.06$&1.00\\
\hline
\rcv& & & & \\
\hline
& LBFGS&$0.6\pm0.3$&$0.47\pm0.08$&&1.24\\
& NCG&$5\pm2$&$0.44\pm0.07$&&10.48\\
& GD&$0.5\pm0.2$&$0.5\pm0.1$&&1.06\\
& LBFGS-P&$0.8\pm0.2$&$0.45\pm0.07$&$0.29\pm0.06$&1.05\\
& NCG-P&$0.8\pm0.2$&$0.46\pm0.07$&$0.30\pm0.06$&1.08\\
& GD-P&$0.8\pm0.1$&$0.45\pm0.07$&$0.29\pm0.06$&1.00\\
\hline
\urldata& & & & \\
\hline
& LBFGS&$6\pm2$&$4.5\pm0.2$&&1.14\\
& NCG&$56\pm15$&$4.6\pm0.3$&&11.10\\
& LBFGS-P&$6.1\pm0.3$&$4.3\pm0.2$&$1.0\pm0.1$&1.00\\
& NCG-P&$5.9\pm0.3$&$4.4\pm0.2$&$1.0\pm0.1$&1.01\\
\hline
\kdda& & & & \\
\hline
& LBFGS&$24\pm5$&$19\pm1$&&1.05\\
& NCG&$30\pm18$&$19\pm1$&&1.38\\
& LBFGS-P&$25\pm2$&$18\pm1$&$4.6\pm0.4$&1.03\\
& NCG-P&$24\pm2$&$18\pm1$&$4.7\pm0.6$&1.03\\
\hline
\kddb& & & & \\
\hline
& LBFGS&$37\pm8$&$28\pm2$&&1.07\\
& NCG&$45\pm28$&$28\pm2$&&1.41\\
& LBFGS-P&$38\pm2$&$27\pm2$&$6.7\pm0.6$&1.00\\
& NCG-P&$36\pm3$&$27\pm2$&$7.1\pm0.8$&1.00\\
\hline
\hline
\end{tabular}

	\end{table}

The advantage of PELS\ for accurate large-scale distributed
line searches is apparent when comparing $\tau_{\ell\text{s}}$ to
$\tau_{\text{fg}}$ in Tab.\;\ref{tab:all}, as well as the difference in
$n_e$ between PELS\ and WA\ algorithms.
For all datasets, $\tau_{\ell\text{s}} < \tau_{\text{fg}}$, however when
the problem size is large, such as for the {\sc kdd-a}\ and {\sc kdd-b}\
datasets, computing the PELS\ coefficients took only a quarter of the
time required to compute $\mathbf{\nabla}\mathcal{L} $.
For all problems, LBFGS-P\ and NCG-P required a lower number
of $n_e$ evaluations in the line search than LBFGS\
and NCG, respectively, which entailed that the average 
$\tau_{\text{tot}}$ was approximately the same for both PELS\ and
WA\ implementations, despite the fact that 
additional work was performed in computing the Taylor coefficients.
In addition, note that $n_e \approx 1$ for all PELS\
algorithms on the datasets considered, which is
particularly effective when contrasted with NCG on the {\sc e}psilon,
{\sc rcv\footnotesize{1}}, and {\sc url}\ datasets: on these problems, NCG
produced poorly scaled search directions requiring many 
line search iterations.
While preconditioning in the NCG algorithm can be used to mitigate the poor
scaling, it requires additional matrix-vector
operations in each iteration, often constructing an
approximation to the Hessian with LBFGS-type updates
\cite{hager2006survey}; we thus find it notable that NCG-P
often had better performance than LBFGS, \emph{without}
the need for preconditioning.
In contrast to NCG, LBFGS\ generally had $n_e \approx
1$ in Tab.\;\ref{tab:all}; as such, the performance gains of LBFGS-P
over LBFGS\ stem principally from reducing the total number of
required iterations by computing more accurate minima in
each line search invocation.

\section{Conclusion.}
\label{sec:conclusion}
In this paper, we have presented the Polynomial Expansion Line Search\ method for
large-scale batch and minibatch optimization
algorithms, applicable to
smooth
loss functions with $L_2$-regularization such as least squares regression, 
logistic regression, and low-rank matrix factorization.
The PELS\ method constructs a truncated Taylor
polynomial expansion of the loss function that may be
minimized quickly and accurately in a neighbourhood of the
expansion point, and additionally has 
coefficients that may be evaluated in parallel with
little communication overhead.
Performance tests with our implementations of
LBFGS, NCG, and GD with PELS\ in the
Apache Spark framework were conducted with a logistic regression model on large
classification datasets on a 16 node cluster with 256
processing cores.
It was found, perhaps surprisingly, that NCG with
PELS often exhibited better convergence and faster
performance than LBFGS\ with a standard Wolfe approximate
line search.
For large datasets, the PELS\ technique also significantly
reduced the number of iterations and time required by the LBFGS\
algorithm
to reach high training
accuracies by factors of 1.8--2.
The PELS technique may be used accelerate
parallel large-scale regression and matrix factorization computations, 
and is applicable to important classes of
smooth optimization problems.
All computer code for this paper is available
through a \texttt{github} repository\footnote{\url{https://github.com/mbhynes}}.

{
	\sloppy
	
		{
								\balance
				\bibliographystyle{ieeetr}
		\bibliography{main} 
	}
 
}
\end{document}